\documentclass[12pt,reqno]{amsart}

\newtheorem*{Theorem}{Theorem}
\newtheorem*{Conjecture}{Conjecture}

\def\({\left(}
\def\){\right)}
\def\[{\left[}
\def\]{\right]}
\def\si{\sigma}
\def\om{\omega}
\def\sgn{\operatorname{sgn}}

\begin{document}

\newbox\Adr
\setbox\Adr\vbox{
\centerline{Institut f\"ur Mathematik der Universit\"at Wien,}
\centerline{Strudlhofgasse 4, A-1090 Wien, Austria.}
\centerline{E-mail: {\tt\footnotesize kratt@ap.univie.ac.at}}
\centerline{WWW: \footnotesize\tt http://www.mat.univie.ac.at/People/kratt}
}

\title[Summation formula for an  $\tilde A_n$ basic hypergeometric series]{Proof of a
summation formula for an  $\tilde A_n$ basic hypergeometric series
conjectured by Warnaar}
\author[C.~Krattenthaler]{C.~Krattenthaler$^\dagger$\\[18pt]\box\Adr}

\address{Institut f\"ur Mathematik der Universit\"at Wien,
Strudlhofgasse 4, A-1090 Wien, Austria.\newline
e-mail: KRATT@Pap.Univie.Ac.At\\
WWW: \tt http://radon.mat.univie.ac.at/People/kratt}

\thanks {$^\dagger$ Research partially supported by the Austrian
Science Foundation FWF, grant P12094-MAT}

\begin{abstract}
A proof of an unusual summation formula for a basic hypergeometric
series associated to the affine root system $\tilde A_n$ that was
conjectured by Warnaar is given. It makes use of
Milne's $A_n$ extension of Watson's transformation, Ramanujan's
$_1\psi_1$-summation, and a determinant evaluation of the author.
In addition, a transformation formula between basic hypergeometric
series associated to the affine root systems $\tilde A_n$
respectively $\tilde A_m$, which generalizes at the same time the above
summation formula and an identity due to Gessel and
the author, is proposed as a conjecture.
\end{abstract}

\maketitle

\begin{section}{Introduction, statement of the result, and of the
conjecture}

The purpose of this note is to prove a summation formula
for a basic hypergeometric
series associated to the affine root system $\tilde A_{n-1}$ that was
conjectured by Warnaar (private communication). (Another frequently
used term for such series is `basic hypergeometric series in
$SU(n)$.' We follow however the terminology for multiple basic
hypergeometric series associated to root systems
as laid down in \cite[Sec.~7]{GeKrAA} and  \cite[Sec.~1]{BhScAA}.
For an overview of the state of the art of this theory and of its
relevance we refer the reader to \cite{MilnAL,BhScAA,VDieAE,KrScAA} and the
references cited therein.)

\begin{Theorem} \label{T1}
Let $n$ be a positive integer, let
$M_1$ and $M_2$ be nonnegative integers, and let $S$ be an
integer with $-M_1\le S\le M_2$. Then
\begin{multline} \label{eq:0}
\sum _{k_1+\dots+k_n=S} ^{}(-1)^{(n-1)S}
q^{\binom {n+1}2\sum _{i=1} ^{n}k_i^2+\sum
_{i=1} ^{n}ik_i}\prod _{1\le i<j\le n} ^{}(1-q^{nk_j-nk_i+j-i})\\
\cdot \prod _{i=1} ^{n}\frac {(q;q)_{M_1+M_2+i-1}}
{(q;q)_{M_1+nk_i+i-1}\,(q;q)_{M_2-nk_i+n-i}}=
q^{(n+1)\binom {S+1}2}\frac {(q;q)_{M_1+M_2}}
{(q;q)_{M_1+S}\,(q;q)_{M_2-S}},
\end{multline}
where, as usual, the shifted $q$-factorial $(a;q)_n$ is defined by
$(a;q)_k:=(1-a)(1-aq)\cdots(1-aq^{k-1})$ if $k>0$, $(a;q)_0:=1$, and
$(a;q)_k:=1/(1-a/q)(1-a/q^2)\cdots(1-aq^{k})$ if $k<0$.
\end{Theorem}

This identity is remarkable, because it essentially\footnote{In fact,
Milne's identity is the $M_2\to\infty$, $M_1=0$ case of \eqref{eq:0}.
However, it is shown in \cite[paragraph before Theorem~22]{GeKrAA},
that, by what is called there the ``rotation trick",
Milne's identity does also imply the $M_2\to\infty$ case of
\eqref{eq:0} (i.e., with $M_1$ arbitrary).
The rotation trick will also be used in our proof of
the Theorem.} reduces to an identity
originally due to Milne \cite[Theorem~1.9]{MilnAG} if we let $M_2$ tend to
infinity.
The proof of Milne's identity in \cite{MilnAG} uses a great
deal of machinery
(in fact a large part of his paper \cite{MilnAG} is devoted
to the proof of this identity), which, apparently, does not allow any
generalization or extension. On the other hand, an elementary,
combinatorial proof of Milne's identity has been given in
\cite[Theorem~22]{GeKrAA}. But, again, it seems impossible to extend
this combinatorial approach to a proof of the above Theorem.

I will prove the above Theorem by an unusual combination of, on the one
hand, classical and, on the other hand, more recent results in classical
analysis. The proof will require
Milne's $A_n$ extension of Watson's transformation
\cite[Theorem~6.1]{MiLiAB}, Ramanujan's classical
$_1\psi_1$-summation (see e.g\@. \cite[Eq.~(5.2.1); Appendix
(II.29)]{GaRaAA}), and a determinant evaluation of the author
\cite[Lemma~2.2]{KratAM}
which is ubiquitous in classcial and combinatorial analysis (cf\@.
\cite[Theorem~26 and the subsequent paragraphs]{KratBN}
for a list of occurrences).

An independent proof of the above Theorem results from an identity
for supernomial coefficients due to Schilling
and Shimozono \cite[Eq.~(6.6)]{ScShAB} (cf\@. \cite[remarks preceding
Eq.~(6.6)]{WarnAI}). I believe that the proof of this paper is
still of interest, because variations of this approach will certainly
turn out to be useful in other cases as well.

A test candidate for the above judgement may be the following
conjectural generalization of the Theorem.
Before I state it precisely, let me
recall that in \cite[Theorem~26]{GeKrAA}
it is shown that Milne's identity (i.e.,
the $M_2\to\infty$ case of the above Theorem) is in fact
part of an infinite {\it hierarchy} of {\it transformation} formulas between
multiple basic hypergeometric {\it of different dimension}. (Such
transformations are, up to now, very rare. Except for Section~8 of
\cite{GeKrAA}, the only occurrence of such transformations that I am
aware of is \cite{KajiAA}.) Since Milne's identity admits the
generalization stated in the above Theorem, an immediate question is
whether or not it is possible to also introduce an additional
parameter into this infinite hierarchy of transformation formulas. On
the basis of computer experiments, there is overwhelming evidence
that this is indeed the case. We state the
formula in the Conjecture below.

\begin{Conjecture} \label{C1}
Let $n$ and $m$ be positive integers, let
$M_1$ and $M_2$ be nonnegative integers, and let $S_1$ and $S_2$ be
integers with $-M_1\le S_1\le M_2$ and $-M_1\le S_2\le M_2$. Then
\begin{multline} \notag
\sum _{k_1+\dots+k_n=S_1} ^{}(-1)^{(n-1)S_1}
q^{\frac {n(n+m)}2 \sum _{i=1} ^{n}k_i^2+m \sum _{i=1} ^{n}i k_i
-m\binom {S_1+1}2-nS_1(S_1+m)/2}\\
\cdot
  \prod _{1\le i<j\le n} ^{}(1-q^{n k_j-n k_i+j-i})
  \prod _{i=1} ^{n}\frac {(q;q)_{M_1+M_2+i-1}} { (q;q)_{M_1-S_1+n k_i+i-1}\,
          (q;q)_{M_2+S_1-n k_i+n-i}}
\end{multline}
\begin{multline} \label{e2}
=\sum _{l_1+\dots+l_m=S_2} ^{}(-1)^{(m-1)S_2}
q^{\frac {m(m+n)}2 \sum _{i=1} ^{m}l_i^2+n \sum _{i=1} ^{m}i l_i
-n\binom {S_2+1}2-mS_2(S_2+n)/2}\kern2cm\\
\cdot
  \prod _{1\le i<j\le m} ^{}(1-q^{m l_j-m l_i+j-i})
  \prod _{i=1} ^{m}\frac {(q;q)_{M_1+M_2+i-1}} { (q;q)_{M_1-S_2+m l_i+i-1}\,
          (q;q)_{M_2+S_2-m l_i+m-i}}.
\end{multline}
\end{Conjecture}

Clearly, our Theorem is the $m=1$ case of this conjecture. Even more
evidence in favour of the conjecture comes from the fact that for
$M_2\to\infty$ it reduces to Theorem~26 in \cite{GeKrAA}.

\medskip
By means of the ``rotation trick" (see \cite[paragraph before
Theorem~22]{GeKrAA} and the first paragraph of the next section), it
can be seen that it suffices to prove the Conjecture for $S_1=S_2=0$.
However, in contrast to our proof of the Theorem,
for a proof of the Conjecture it will not be sufficient to apply Milne's
$A_{n-1}$ extension of Watson's transformation. Perhaps one has to start
with a higher order transformation formula,
for example, with one of the $A_{n-1}$ extensions of Bailey's
very-well-poised $_{10}\phi_9$-transformation formula from
\cite{MiNeAA}.

\end{section}

\begin{section} {Proof of the Theorem}

First of all, analogously to the remark of the last paragraph of the
previous section, I claim that it is enough to prove \eqref{eq:0} for
$S=0$, i.e.,
\begin{multline} \label{eq:1}
\sum _{k_1+\dots+k_n=0} ^{}q^{\binom {n+1}2\sum _{i=1} ^{n}k_i^2+\sum
_{i=1} ^{n}ik_i}\prod _{1\le i<j\le n} ^{}(1-q^{nk_j-nk_i+j-i})\\
\cdot \prod _{i=1} ^{n}\frac {(q;q)_{M_1+M_2+i-1}}
{(q;q)_{M_1+nk_i+i-1}\,(q;q)_{M_2-nk_i+n-i}}=\frac {(q;q)_{M_1+M_2}}
{(q;q)_{M_1}\,(q;q)_{M_2}}.
\end{multline}
This is seen by resorting to the ``rotation trick" \cite[paragraph before
Theorem~22]{GeKrAA}. Let us assume that we already proved
\eqref{eq:1}. Let $S$ be some fixed
integer. Division of $S$ by $n$ gives a unique
representation $S=Qn+R$ where $Q,R$ are integers with $0\le R<
n$. Then in \eqref{eq:1} replace $k_1$ by $k_{1+R}-Q$, \dots, $k_{n-R}$ by
$k_n-Q$, $k_{n-R+1}$ by $k_1-Q-1$, \dots, $k_n$ by $k_R-Q-1$. So the
effect is a rotation of the summation indices, combined with a certain
shift. If we rewrite \eqref{eq:1} after these replacements and finally replace
$M_1$ by $M_1+S$ and $M_2$ by $M_2-S$,
we obtain \eqref{eq:0} after some simplification.

Next, I claim that it is enough to prove \eqref{eq:1} for $M_1\equiv 0$ mod
$n$. To see this, suppose that $M_2$ is given. Multiply both sides of
\eqref{eq:1} by $\prod _{i=1} ^{n}(q^{M_1+M_2+i};q)_n$ and write the
result in the form
\begin{multline} \label{eq:2}
\sum _{k_1+\dots+k_n=0} ^{}q^{\binom {n+1}2\sum _{i=1} ^{n}k_i^2+\sum
_{i=1} ^{n}ik_i}\prod _{1\le i<j\le n} ^{}(1-q^{nk_j-nk_i+j-i})\\
\times \prod _{i=1} ^{n}\frac {(q^{M_1+nk_i+i};q)_{M_2-nk_i+n}}
{(q;q)_{M_2-nk_i+n-i}}=
\frac {(q^{M_1+1};q)_{M_2}} {(q;q)_{M_2}}
\prod _{i=1} ^{n}(q^{M_1+M_2+i};q)_n.
\end{multline}
Both sides are most obviously polynomials in $q^{M_1}$, of degree at
most $n^2(n+M_2)$, because, in the summation, each $k_i$ is bounded above
by $1+M_2/n$, and, hence, bounded below by $-(n-1)(1+M_2/n)$. A
polynomial is uniquely determined by its evaluation at enough points,
certainly at infinitely many points. Therefore, if \eqref{eq:2} is
true for all $M_1\equiv 0$ mod $n$ then it is true for all $M_1$.
Since \eqref{eq:2} and \eqref{eq:1} are equivalent, the same applies
to \eqref{eq:1}.

Now, choose some $M_1\equiv0$ mod $n$. If we want to prove
\eqref{eq:1} for this particular $M_1$, then an analogous argument
shows that it is enough to prove it for all $M_2\equiv0$ mod $n$.

Summarizing, it is sufficient to prove \eqref{eq:1} for
$M_1\equiv M_2\equiv0$ mod $n$. Therefore, for the rest of the proof, we
assume that this congruence condition is satisfied.

To begin with, let us rewrite the left-hand side of \eqref{eq:1}
by replacing $k_i$ by
$k_i-M_1/n$, $i=1,2,\dots,n$, and performing some rearrangement of terms,
\begin{multline} \label{eq:3}
(-1)^{nM_1}q^{M_1(M_1n-M_1+2nM_2+2n^2-1)/2}\\
\times
\sum _{k_1+\dots+k_n=M_1} ^{}q^{{\frac {n} {2}\sum _{i=1}
^{n}k_i^2-(n-1)\sum _{i=1} ^{n}ik_i}}
\prod _{1\le i<j\le n} ^{}\frac {(1-q^{nk_j-nk_i+j-i})} {1-q^{j-i}}
\prod _{i=1} ^{n}\frac {(q^{-M_1-M_2-n+i};q)_{nk_i}}
{(q^i;q)_{nk_i}}.
\end{multline}

Next we want to apply a limiting case of Milne's $A_n$
Watson transformation \cite[Theorem~6.1]{MiLiAB},
\begin{align} \notag
&\sum _{k_1,\dots,k_l\ge0} ^{}\bigg(\prod _{1\le r<s\le l} ^{}\frac {1-\frac
{x_r} {x_s}q^{k_r-k_s}} {1-\frac {x_r} {x_s}}\bigg)
\bigg(\prod _{i=1} ^{l}\frac
{1-\frac {x_i} {x_l}aq^{k_i+(k_1+\dots+k_l)}} {1-\frac {x_i}
{x_l}a}\bigg)
\bigg(\prod _{r=1} ^{l}\prod _{s=1} ^{l}\frac {(\frac {x_r}
{x_s}q^{-N_s};q)_{k_r}} {(q\frac {x_r} {x_s};q)_{k_r}}\bigg)
\\
\notag
&\quad \times
\bigg(\prod _{i=1} ^{l}\frac {(\frac {x_i}
{x_l}a;q)_{k_1+\dots+k_l}} {(\frac {x_i}
{x_l}aq^{1+N_i};q)_{k_1+\dots+k_l}}\bigg)
\bigg(\prod _{i=1} ^{l}\frac {(\frac {x_i}
{x_l}c;q)_{k_i}\,(\frac {x_i}
{x_l}d;q)_{k_i}} {(\frac {x_i} {x_l}\frac {aq} {b};q)_{k_i}\,
(\frac {x_i} {x_l}\frac {aq} {e};q)_{k_i}}\bigg)
\\
\notag
&\quad  \times
\frac {(b;q)_{k_1+\dots+k_l}\,(e;q)_{k_1+\dots+k_l}} {(\frac {aq}
{c};q)_{k_1+\dots+k_l}\,(\frac {aq} {d};q)_{k_1+\dots+k_l}}
\(\frac {a^2q^{1+N_1+\dots+N_l}} {bcde}\)^{k_1+\dots+k_l}
q^{\sum _{i=1} ^{l}ik_i}\\
\notag
&=
\frac {(aq/de;q)_{N_1+\dots+N_l}} {(aq/d;q)_{N_1+\dots+N_l}}
\bigg(\prod _{i=1} ^{l}\frac {(\frac {x_i} {x_l}aq;q)_{N_i}}
{(\frac {x_i} {x_l}aq/e;q)_{N_i}}\bigg)
\sum _{k_1,\dots,k_l\ge0} ^{}
q^{\sum _{i=1} ^{l}ik_i}
\bigg(\prod _{1\le r<s\le l} ^{}\frac {1-\frac
{x_r} {x_s}q^{k_r-k_s}} {1-\frac {x_r} {x_s}}\bigg)\\
\label{eq:MiLiAB}
&\quad \times
\bigg(\prod _{r=1} ^{l}\prod _{s=1} ^{l}\frac {(\frac {x_r}
{x_s}q^{-N_s};q)_{k_r}} {(q\frac {x_r} {x_s};q)_{k_r}}\bigg)
\bigg(\prod _{i=1} ^{l}\frac {(\frac {x_i}
{x_l}d;q)_{k_i}} {(\frac {x_i}
{x_l}\frac {aq} {b};q)_{k_i}}\bigg)
\frac {(\frac {aq} {bc};q)_{k_1+\dots+k_l}\,(e;q)_{k_1+\dots+k_l}}
{(\frac {aq}
{c};q)_{k_1+\dots+k_l}\,(\frac {de}
{a}q^{-N_1-\dots-N_l};q)_{k_1+\dots+k_l}},
\end{align}
where $N_1,\dots,N_l$ are nonnegative integers.
For convenience, let us set $k_{l+1}=M-k_1-\dots-k_l$ and $a=x_l/q^Mx_{l+1}$,
so that \eqref{eq:MiLiAB} becomes
\begin{align} \notag
&q^M\bigg(\prod _{i=1} ^{l}\frac {(q\frac {x_{l+1}} {x_i};q)_M}
{(\frac {x_{l+1}} {x_i}q^{-N_i};q)_M}\bigg)
\frac {(q^{1-M}/b;q)_M\,(q^{1-M}/e;q)_M} {(\frac {x_{l+1}}
{x_l}c;q)_M\,(\frac {x_{l+1}} {x_l}d;q)_M}
\prod _{i=1} ^{l}\frac {(1-\frac {x_i} {x_{l+1}})}
{(1-\frac {x_i} {q^Mx_{l+1}})}
\\
\notag
&\quad \times
\sum _{k_1+\dots+k_l+k_{l+1}=M} ^{}
q^{-\sum _{i=1} ^{l+1}ik_i}
\bigg(\prod _{1\le r<s\le l+1} ^{}\frac {1-\frac
{x_s} {x_r}q^{k_s-k_r}} {1-\frac {x_s} {x_r}}\bigg)
\bigg(\prod _{r=1} ^{l+1}\prod _{s=1} ^{l}\frac {(\frac {x_r}
{x_s}q^{-N_s};q)_{k_r}} {(q\frac {x_r} {x_s};q)_{k_r}}\bigg)
\\
\notag
&\kern4cm \times
\bigg(\prod _{i=1} ^{l+1}\frac {(\frac {x_i}
{x_l}c;q)_{k_i}\,(\frac {x_i}
{x_l}d;q)_{k_i}} {(\frac {x_i} {x_{l+1}}\frac {q^{1-M}} {b};q)_{k_i}\,
(\frac {x_i} {x_{l+1}}\frac {q^{1-M}} {e};q)_{k_i}}\bigg)
\end{align}
\begin{align}
\notag
&=
\frac {(q^{1-M}x_l/x_{l+1}de;q)_{N_1+\dots+N_l}}
{(q^{1-M}x_l/x_{l+1}d;q)_{N_1+\dots+N_l}}
\bigg(\prod _{i=1} ^{l}\frac {(\frac {x_i} {x_{l+1}}q^{1-M};q)_{N_i}}
{(\frac {x_i} {x_{l+1}}q^{1-M}/e;q)_{N_i}}\bigg)\\
\notag
&\quad \times
\sum _{k_1,\dots,k_l\ge0} ^{}
q^{\sum _{i=1} ^{l}ik_i}
\bigg(\prod _{1\le r<s\le l} ^{}\frac {1-\frac
{x_r} {x_s}q^{k_r-k_s}} {1-\frac {x_r} {x_s}}\bigg)
\bigg(\prod _{r=1} ^{l}\prod _{s=1} ^{l}\frac {(\frac {x_r}
{x_s}q^{-N_s};q)_{k_r}} {(q\frac {x_r} {x_s};q)_{k_r}}\bigg)
\\
\label{eq:MiLiAB-mod}
&\kern2cm \times
\bigg(\prod _{i=1} ^{l}\frac {(\frac {x_i}
{x_l}d;q)_{k_i}} {(\frac {x_i}
{x_{l+1}}\frac {q^{1-M}} {b};q)_{k_i}}\bigg)
\frac {(\frac {q^{1-M}x_l} {x_{l+1}bc};q)_{k_1+\dots+k_l}\,(e;q)_{k_1+\dots+k_l}}
{(\frac {q^{1-M}x_l}
{x_{l+1}c};q)_{k_1+\dots+k_l}\,
(\frac {deq^Mx_{l+1}} {x_l}q^{-N_1-\dots-N_l};q)_{k_1+\dots+k_l}}.
\end{align}
In this identity we replace $q$ by $q^n$.
Then we set $l=n-1$, $M=M_1$, $x_i=q^i$ for $i=1,2,\dots,l+1$,
$b=q^{-nM_1}$, $d=\delta q^{-M_1-M_2-1}$, $N_i=(M_1+M_2)/n$. Next we
multiply both sides by $(1-\delta)$ (this cancels one factor in the term
$(x_{l+1}d/x_{l};q)_M\sim (\delta q^{-M_1-M_2};q^n)_{M_1}$ in the denominator
of the left-hand side of \eqref{eq:MiLiAB-mod} and one factor in the term
$(q^{1-M}x_l/x_{l+1}d;q)_{N_1+\dots+N_l}\sim
(q^{n-nM_1+M_1+M_2}/\delta;q^n)_{(M_1+M_2)/n}$ in the denominator
of the right-hand side of \eqref{eq:MiLiAB-mod}).
Finally, we let $\delta\to1$,
$c\to\infty$, and $e\to\infty$. This reduces \eqref{eq:MiLiAB-mod} to
the following transformation formula,
\begin{align} \notag
&q^{nM_1(n-M_1/2-1)}
\frac {(q;q)_{nM_1}}
{(q^{-M_1-M_2};q)_{M_1+M_2}\,(q;q)_{nM_1-M_1-M_2-1}}\\
\notag
&\quad \times
\sum _{k_1+\dots+k_n=M_1} ^{}
q^{\frac {n} {2}\sum _{i=1} ^{n}k_i^2-(n-1)\sum _{i=1} ^{n}ik_i}
\bigg(\prod _{1\le r<s\le n} ^{}\frac {1-q^{nk_s-nk_r+s-r}}
{1-q^{s-r}}\bigg)
\\
\notag
&\kern4cm \times
\bigg(\prod _{r=1} ^{n}\frac {(q^{-M_1-M_2-n+r};q)_{nk_r}}
{(q^r;q)_{nk_r}}\bigg)
\\
\notag
&=-\frac {(q^{-nM_1};q)_{M_1+M_2+1}}
{(q^{-nM_1};q^n)_{M_1}\,(q^n;q^n)_{M_2}}
\sum _{k_1,\dots,k_{n-1}\ge0} ^{}
\bigg(\prod _{1\le r<s\le n-1} ^{}\frac {1-q^{nk_r-nk_s+r-s}}
{1-q^{r-s}}\bigg)\\
&\kern2cm \times
\bigg(\prod _{r=1} ^{n-1}\frac {(q^{-M_1-M_2-n+r};q)_{nk_r}}
{(q^r;q)_{nk_r}}\(q^{M_2}\)^{nk_r}\bigg)
q^{n\sum _{i=1} ^{n-1}ik_i}.
\label{eq:MiLiAB-lim}
\end{align}
The series on the left-hand side of \eqref{eq:MiLiAB-lim} is
exactly the series in \eqref{eq:3}. What the transformation
\eqref{eq:MiLiAB-lim} does with this series is, in some sense which will
become more transparent below, that it ``entangles" the
summation indices. Thus, we obtain the following expression for
the left-hand side of \eqref{eq:1},
\begin{multline}
(-1)^{M_1}q^{-\binom {M_1+1}2}
\frac {(q;q)_{M_1+M_2}}
{(q^{-nM_1};q^n)_{M_1}\,(q^n;q^n)_{M_2}}\\
\times
\sum _{k_1,\dots,k_{n-1}\ge0} ^{}
\bigg(\prod _{1\le r<s\le n-1} ^{}\frac {1-q^{nk_r-nk_s+r-s}}
{1-q^{r-s}}\bigg)\kern3cm\\
\times
\bigg(\prod _{r=1} ^{n-1}\frac {(q^{-M_1-M_2-n+r};q)_{nk_r}}
{(q^r;q)_{nk_r}}\(q^{M_2}\)^{nk_r}\bigg)
q^{n\sum _{i=1} ^{n-1}ik_i}.
\label{eq:4}
\end{multline}
(The sign $(-1)^{M_1}$ is no misprint since our assumption
$M_1\equiv 0$ mod $n$ implies $nM_1\equiv M_1$ mod 2.)

The next task is to split the sum in \eqref{eq:4} into many pieces,
each of which being a product of $n-1$ one-dimensional summations.
This is done by replacing the product over $1\le s<r\le n-1$ by a
Vandermonde determinant. More precisely, we have
\begin{align*}
\prod _{1\le r<s\le n-1} ^{} (1-q^{nk_r-nk_s+r-s})
&=q^{-\sum _{i=1} ^{n-1}(i-1)(nk_i+i)}
\prod _{1\le r<s\le n-1} ^{} (q^{nk_s+s}-q^{nk_r+r}) \\
&=q^{-\sum _{i=1} ^{n-1}(i-1)(nk_i+i)}
\det_{1\le i,j\le n-1}\(\(q^{nk_i+i}\)^{j-1}\)\\
&=q^{-n\sum _{i=1} ^{n-1}ik_i+n\sum _{i=1} ^{n-1}k_i-2\binom n3}
\sum _{\si\in S_{n-1}} ^{}\sgn \si \prod _{i=1}
^{n-1}q^{(\si(i)-1)(nk_i+i)}.
\end{align*}
Hence, the sum in \eqref{eq:4} equals
\begin{multline}
(-1)^{\binom {n-1}2}\bigg(\prod _{i=1} ^{n-1}\frac {1} {(q;q)_{i-1}}\bigg)
\sum _{\si\in S_{n-1}} ^{}\sgn \si\, q^{-\binom n3}q^{\sum _{i=1}
^{n-1}i(\si(i)-1)}\\
\times
\prod _{i=1} ^{n-1}
\bigg(\sum _{k_i\ge0} ^{}
\frac {(q^{-M_1-M_2-n+i};q)_{nk_i}}
{(q^i;q)_{nk_i}}\(q^{M_2+\si(i)}\)^{nk_i}\bigg).
\label{eq:5}
\end{multline}
The next ingredient is Ramanujan's $_1\psi_1$-summation (see
\cite[(5.2.1)]{GaRaAA}),
\begin{equation} \label{eq:1psi1}
\sum _{k=-\infty} ^{\infty}\frac {(a;q)_k} {(b;q)_k}z^k=
\frac {(q;q)_\infty\,(b/a;q)_\infty\,(az;q)_\infty\,(q/az;q)_\infty}
{(b;q)_\infty\,(q/a;q)_\infty\,(z;q)_\infty\,(b/az;q)_\infty}.
\end{equation}
Each of the inner sums in \eqref{eq:5} is an $n$-section of a special case
of the left-hand side of \eqref{eq:1psi1}. (To be precise, it is the
special case $a=q^{-M_1-M_2-n+i}$, $b=q^i$, and $z=q^{M_2+\si(i)}$.)
Thus, \eqref{eq:5} simplifies to
\begin{multline}
(-1)^{\binom {n-1}2}\bigg(\prod _{i=1} ^{n-1}\frac {1} {(q;q)_{i-1}}\bigg)
\sum _{\si\in S_{n-1}} ^{}\sgn \si\, q^{-\binom n3}q^{\sum _{i=1}
^{n-1}i(\si(i)-1)}\\
\times
\prod _{i=1} ^{n-1}
\Bigg(\frac {1} {n}\sum _{\ell_i=0} ^{n-1}\frac
{(q;q)_\infty\,(q^{M_1+M_2+n};q)_\infty\,
(q^{i+\si(i)-M_1-n}\om^{\ell_i};q)_\infty\,
(q^{1-i-\si(i)+M_1+n}\om^{-\ell_i};q)_\infty}
{(q^i;q)_\infty\,(q^{1-i+M_1+M_2+n};q)_\infty\,
(q^{M_2+\si(i)}\om^{\ell_i};q)_\infty\,
(q^{-\si(i)+M_1+n}\om^{-\ell_i};q)_\infty}\Bigg),
\label{eq:6}
\end{multline}
where $\om$ denotes a primitive $n$-th root of unity. An immediate
observation is that if any $\ell_i$ equals 0 then the corresponding
summand vanishes, because of the term
$$(q^{i+\si(i)-M_1-n}\om^{\ell_i};q)_\infty\,
(q^{1-i-\si(i)+M_1+n}\om^{-\ell_i};q)_\infty$$
in the numerator. Hence,
we may as well sum over $\ell_i$ from $1$ to $n-1$,
$i=1,2,\dots,n-1$.

Some manipulation transforms \eqref{eq:6} into
\begin{multline}
(-1)^{\binom {n-1}2}\frac {1} {n^{n-1}}
\prod _{i=1} ^{n-1}\frac {1} {(q^{1-i+M_1+M_2+n};q)_{i-1}}
\sum _{\si\in S_{n-1}} ^{}\sgn \si\, \\
\cdot\Bigg(\sum _{\ell_1,\dots,\ell_{n-1}=1} ^{n-1}
\bigg(\prod _{i=1} ^{n-1}\frac {(q^{-M_1}\om^{\ell_i};q)_{\infty}}
{(q^{M_2+1}\om^{\ell_i};q)_\infty}\om^{\ell_i(n-i-\si(i))}\bigg)\\
\cdot
(q^{M_2+1}\om^{\ell_i};q)_{\si(i)-1}\,
(q^{M_1+1}\om^{-\ell_i};q)_{n-\si(i)-1}
\Bigg).
\label{eq:7}
\end{multline}

Now it is not difficult to see that if $\ell_{r}=\ell_s$, $r\ne s$,
then the summand corresponding to the permutation $\si$ cancels with
the summand corresponding to the permutation $\si\circ (rs)$. (Here,
$(rs)$ denotes the transposition which interchanges $r$ and $s$.)
Therefore the only summands which survive this cancellation are those
where the summation indices $\ell_1,\ell_2,\dots,\ell_{n-1}$ are a
permutation of $\{1,2,\dots,n-1\}$. Thus, \eqref{eq:7} reduces to
\begin{multline}
(-1)^{\binom {n-1}2}\frac {1} {n^{n-1}}
\frac {(q^{-nM_1};q^n)_{M_1}\, (q^n;q^n)_\infty\,
(q^{M_2+1};q)_\infty} {(q^{-M_1};q)_{M_1}\, (q;q)_\infty\,
(q^{nM_2+n};q^n)_\infty}
\prod _{i=1} ^{n-1}\frac {(1-\om^i)} {(q^{1-i+M_1+M_2+n};q)_{i-1}}\\
\times
\sum _{\si,\tau\in S_{n-1}} ^{}(\sgn \si)\,
\om^{\tau(i)(n-i-\si(i))}
(q^{M_2+1}\om^{\tau(i)};q)_{\si(i)-1}\,
(q^{M_1+1}\om^{-\tau(i)};q)_{n-\si(i)-1}
\\
=
(-1)^{\binom {n-1}2}\frac {1} {n^{n-1}}
\frac {(q^{-nM_1};q^n)_{M_1}\, (q^n;q^n)_\infty\,
(q^{M_2+1};q)_\infty} {(q^{-M_1};q)_{M_1}\, (q;q)_\infty\,
(q^{nM_2+n};q^n)_\infty}
\prod _{i=1} ^{n-1}\frac {(1-\om^i)} {(q^{1-i+M_1+M_2+n};q)_{i-1}}\\
\times
\sum _{\tau\in S_{n-1}} ^{}
\om^{\tau(i)(n-i)}\det_{1\le i,j\le n-1}\big(\om^{-j\tau(i)}
(q^{M_2+1}\om^{\tau(i)};q)_{j-1}\,
(q^{M_1+1}\om^{-\tau(i)};q)_{n-j-1}\big).
\label{eq:8}
\end{multline}

The determinant is easily evaluated with the help of the
determinant lemma \cite[Lemma~2.2]{KratAM},
\begin{multline}
\label{lem:Kratt}
\det_{1\le i,j\le n}\Big((X_i+A_n)\cdots(X_i+A_{j+1})
(X_i+B_j)\cdots (X_i+B_2)\Big)\\
\hskip2cm =\prod _{1\le i<j\le n} ^{}(X_i-X_j)\prod _{2\le i\le j\le n}
^{}(B_i-A_j),
\end{multline}
where $X_1,\dots,X_n$, $A_2,\dots,A_n$, and $B_2,
\dots B_n$ are arbitrary indeterminates.
In order to apply \eqref{lem:Kratt}, we rewrite the determinant in
\eqref{eq:8} as
\begin{multline*}
\det_{1\le i,j\le n-1}\big(\om^{-j\tau(i)}
(q^{M_2+1}\om^{\tau(i)};q)_{j-1}\,
(q^{M_1+1}\om^{-\tau(i)};q)_{n-j-1}\big)\\
=(-1)^{\binom {n-1}2}\prod _{i=1}
^{n-1}\om^{-\tau(i)}q^{(M_1+1)+(M_1+2)+\dots+(M_1+n-i-1)}
\kern4cm\\
\det_{1\le i,j\le n-1}\big((\om^{-\tau(i)}-q^{M_2+1})
(\om^{-\tau(i)}-q^{M_2+2})\dots(\om^{-\tau(i)}-q^{M_2+j-1})\\
\cdot (\om^{-\tau(i)}-q^{-M_1-n+j+1})
(\om^{-\tau(i)}-q^{-M_1-n+j+2})\cdots
(\om^{-\tau(i)}-q^{-M_1-1})\big).
\end{multline*}
Now the determinant evaluation \eqref{lem:Kratt} applies with
$X_i=\om^{-\tau(i)}$, $A_j=-q^{-M_1-n+j}$, and $B_j=-q^{M_2+j-1}$.
If the resulting expression is substituted back into \eqref{eq:8}, we obtain
\begin{multline}
\frac {1} {n^{n-1}}
\frac {(q^{-nM_1};q^n)_{M_1}\, (q^n;q^n)_\infty\,
(q^{M_2+1};q)_\infty} {(q^{-M_1};q)_{M_1}\, (q;q)_\infty\,
(q^{nM_2+n};q^n)_\infty}
\prod _{i=1} ^{n-1}{(1-\om^i)} \\
\times
\sum _{\tau\in S_{n-1}} ^{}
\om^{\tau(i)(n-i-1)}
\prod _{1\le i<j\le n-1} ^{}(\om^{-\tau(i)}-\om^{-\tau(j)})\\
=
\frac {1} {n^{n-1}}
\frac {(q^{-nM_1};q^n)_{M_1}\, (q^n;q^n)_\infty\,
(q^{M_2+1};q)_\infty} {(q^{-M_1};q)_{M_1}\, (q;q)_\infty\,
(q^{nM_2+n};q^n)_\infty}
\prod _{i=1} ^{n-1}{(1-\om^i)}
\prod _{1\le i<j\le n-1} ^{}(\om^{-i}-\om^{-j})
\\
\times
\sum _{\tau\in S_{n-1}} ^{}(\sgn\tau)\,
\om^{\tau(i)(n-i-1)}.
\label{eq:9}
\end{multline}
The sum over permutations in the last line is just a Vandermonde
determinant, and as such easily evaluated. If we substitute this in
\eqref{eq:9}, the resulting expression for the sum in \eqref{eq:4}, we obtain
\begin{equation}
\frac {1} {n^{n-1}}
\frac {(q;q)_{M_1+M_2}}
{(q;q)_{M_1}\,(q;q)_{M_2}}
\prod _{i=1} ^{n-1}{(1-\om^i)}
\prod _{1\le i<j\le n-1} ^{}(\om^{-i}-\om^{-j})(\om^{i}-\om^{j})
\label{eq:10}
\end{equation}
for the left-hand side of \eqref{eq:1}. Clearly, there holds
$$\prod _{i=1} ^{n-1}{(1-\om^i)} =n,$$
because it is the limit $\lim_{z\to 1}(1-z^n)/(1-z)$.
Moreover, we have
\begin{align*}
\prod _{1\le i<j\le n-1} ^{}&(\om^{i}-\om^{j})(\om^{-i}-\om^{-j})
=\prod _{i=1} ^{n-1}(1-\om)\cdots (1-\om^{i-1})
\prod _{i=1} ^{n-1}(1-\om^{-1})\cdots (1-\om^{-i+1})\\
&=\prod _{i=1} ^{n-2}(1-\om)\cdots (1-\om^{i})
\prod _{i=1} ^{n-2}(1-\om^{n-1})\cdots (1-\om^{i+1})\\
&=\prod _{i=1} ^{n-2}(1-\om)\cdots (1-\om^{n-1})
=n^{n-2},
\end{align*}
in view of the previous observation.
Thus, \eqref{eq:10} does indeed reduce to the right-hand side of
\eqref{eq:1}. In view of the remarks of the first paragraph of this
section, the proof of the theorem is complete.
\qed
\end{section}

\end{document}